\documentclass[a4paper]{amsart}
\usepackage{eucal}
\usepackage[latin1]{inputenc}

\newtheorem{contar}{contar}
\newtheorem{proposition}[contar]{Proposition}

\newtheorem{theorem}[contar]{Theorem}
\newtheorem{definition}[contar]{Definition}

\newtheorem{remarkth}[contar]{Remark}
\newenvironment{remark}{\begin{remarkth}\upshape}{\hfill$\diamond$\end{remarkth}}

\def\emph#1{{\bfseries\itshape{#1}}}
\newcommand{\qquand}{\qquad\text{and}\qquad}
\newcommand{\quand}{\quad\text{and}\quad}
\newcommand{\spC}{^{\scriptscriptstyle C}}
\newcommand{\XsigmaJ}{{X_\sigma^{\scriptscriptstyle(1)}}}
\def\p#1{\mathaccent19{#1}}
\newcommand{\pd}[2]{\frac{\partial#1}{\partial#2}} 
\newcommand{\dpd}[2]{{\displaystyle\pd{#1}{#2}}}
\newcommand{\at}[1]{\Big\vert_{#1}}
\newcommand{\Real}{\mathbb{R}}
\newcommand{\set}[2]{\left\{\,#1\left.\vphantom{#1#2}\,\right\vert\,#2\,
                \right\}}
\newcommand{\map}[3]{#1\colon#2\rightarrow#3}
\newcommand{\cinfty}[1]{C^{\scriptscriptstyle\infty}(#1)}   
\renewcommand{\sec}[1]{\operatorname{Sec}(#1)} 
\newcommand{\vectorfields}[1]{\mathfrak{X}(#1)}

\newcommand{\id}{{\operatorname{id}}}
\newcommand{\tr}{\operatorname{tr}}
\newcommand{\ext}[2][]{\bigwedge\nolimits^{#1}{#2}}

\newcommand{\pai}[2]{\langle{#1},{#2}\rangle} 
\newcommand{\pb}{^\star} % pullback (algebraic meaning)

\newcommand{\Lprol}[1]{\mathcal{L}#1}
\newcommand{\Jprol}[1]{\mathcal{J}#1}
\newcommand{\prEM}{\tau^E_M}
\newcommand{\prFN}{\tau^F_N}
\newcommand{\prEF}{\pi}
\newcommand{\prMN}{\nu}
\newcommand{\Lpi}[1][]{\mathcal{L}_{#1}\pi}
\newcommand{\Jpi}[1][]{\mathcal{J}_{#1}\pi}
\newcommand{\Vpi}[1][]{\mathcal{V}_{#1}\pi}

\newcommand{\prLM}{\tilde{\pi}_{10}}
\newcommand{\prVM}{\boldsymbol{\pi}_{10}}
\newcommand{\prJM}{\pi_{10}}
\newcommand{\prJN}{\pi_1}
\newcommand{\dLpi}[1][]{\mathcal{L}_{#1}^*\!\pi}

\newcommand{\Mor}[1]{\mathcal{M}(#1)} 
\newcommand{\V}{\mathcal{V}}

\begin{document}

\title{Classical field theory on Lie algebroids: Variational aspects}

\author[E.\ Mart\'{\i}nez]{Eduardo Mart\'{\i}nez}
\address{Eduardo Mart\'{\i}nez:
Departamento de Matem\'atica Aplicada,
Facultad de Ciencias,
Universidad de Zaragoza,
50009 Zaragoza, Spain}
\email{emf@unizar.es}

\thanks{Talk delivered at the 9th International Conference on Differential Geometry and its Applications, Prague, September 2004}

\thanks{Partial financial support from MICYT grant  BFM2003-02532 is acknowledged}


\begin{abstract}
The variational formalism for classical field theories is extended to the setting of Lie algebroids. Given a Lagrangian function we study the problem of finding critical points of the action functional when we restrict the fields to be morphisms of Lie algebroids. In addition to the standard case, our formalism includes as particular examples the case of systems with symmetry (covariant Euler-Poincaré and Lagrange Poincaré cases), Sigma models or Chern-Simons theories. 
\end{abstract}

\maketitle

\section{Introduction}

By using the geometry of Lie algebroids, Weinstein~\cite{Weinstein} showed that it is possible to give a common description of the most interesting classical mechanical systems. Later, this theory was extended to time-dependent Classical Mechanics in~\cite{SaMeMa1,MaMeSa}, by introducing a generalization of the notion of Lie algebroid when the bundle is no longer a vector bundle but an affine bundle. The particular case considered in~\cite{SaMeMa2} is the analog of a first jet bundle of a bundle $M\to\Real$, and hence it is but a particular case of a first order field theory for a 1-dimensional space-time. Therefore, it is natural to investigate whether it is possible to extend our formalism to the case of a general field theory, where the space-time manifold is no longer one dimensional.

Generally speaking, there are three different but closely related aspects in the analysis of first order field theories: the \textsl{variational approach}, which leads to the Euler-Lagrange equations, the \textsl{multisymplectic formalism} on the first jet bundle, and the \textsl{infinite-dimensional approach} on the space of Cauchy data. The aim of this paper is to study several aspects of the variational approach to the theory formulated in the context of Lie algebroids, and we will leave for the future the study of the multisymplectic and the infinite-dimensional formalism.

The standard geometrical approach to the Lagrangian description of first order Classical Field Theories~\cite{CaCrIb,Barna-L,GIMMSY1,BiSnFi,GiMaSa} is based on the canonical structures on the first order jet bundle~\cite{Saunders} of a fiber bundle whose sections are the fields of the theory. Thinking of a Lie algebroid as a substitute for the tangent bundle to a manifold, the analog of the field bundle to be considered here is a surjective morphism of Lie algebroids $\map{\pi}{E}{F}$. In the standard theory, a 1-jet of a section of a bundle is but the tangent map to that section at the given point, and therefore it is a linear map between tangent spaces which has to be a section of the tangent of the projection map. In our theory, the analog object is a linear map from a fiber of $F$ to a fiber of $E$ which is a section of the projection $\pi$. The space $\Jpi$ of these maps has the structure of an affine bundle and it is the space where our theory is developed.

In the case $F=TN$ we can set a variational problem consisting in finding those morphisms of Lie algebroids which are critical points of the integral of a Lagrangian function defined on $\Jpi$. By the choice of some special variations determined by the underlying geometry, we can find the Euler-Lagrange equations.

The situation is to be compared with the so-called covariant reduction theory~\cite{CaRaSh,CaGaRa,CaRa}, where the field equations are obtained by reducing the variational problem, i.e. by restricting the variations to those coming from variations for the original unreduced problem. Therefore, the equations are of different form in each case, Euler-Poincaré, Lagrange-Poincaré for a system with symmetry and (a somehow different) Lagrange-Poincaré for semidirect products. In contrast, our theory includes all these cases as particular cases of a common setting.

We finally mention that related ideas to those explained in this paper have been recently considered by Strobl and coworkers~\cite{Strobl,BoKoSt} in the so-called off-shell theory. In this respect our theory should be considered as the on-shell counterpart for the cases considered there.

\smallskip

The organization of the paper is as follows. In section~\ref{preliminaries} we will recall some differential geometric structures related to Lie algebroids, such as the exterior differential and the flow defined by a section. In section~\ref{jetoids} we define the analog of the manifold of 1-jets in this algebraic setting, which is the space where the Lagrangian of the theory is defined. In section~\ref{morphisms} we find local conditions for a section of our bundle to be a morphism of Lie algebroids. In section~\ref{complete.lift} we define the concept of complete lift of a section, which plays the role of the complete lift or first jet prolongation of a vector field in the standard theory. In section~\ref{variational} we state the constrained variational problem we want to solve and we find the Euler-Lagrange equations. Finally, in section~\ref{examples} we show some examples which illustrate the results in this paper.

%-------------------------------------------------------------------
\subsection*{Notation}
All manifolds and maps are taken in the $C^\infty$ category. The set of smooth functions on a manifold $M$ will be denoted by $\cinfty{M}$. The set of smooth vectorfields on a manifold $M$ will be denoted $\vectorfields{M}$. The set of smooth sections of a fiber bundle $\map{\pi}{B}{M}$  will be denoted $\sec{\pi}$ or $\sec{B}$ when there is no possible confusion. 

The tensor product of a vector bundle $\map{\tau}{E}{M}$ by itself $p$ times will be denoted by $\map{\tau^{\otimes p}}{E^{\otimes p}}{M}$, and similarly the exterior power will be denoted by  $\map{\tau^{\wedge p}}{E^{\wedge p}}{M}$. The set of sections of the dual $\map{\tau^*{}^{\wedge p}}{(E^*)^{\wedge p}}{M}$ will be denoted by $\ext[p]{E}$. For $p=0$ we have $\ext[0]{E}=\cinfty{M}$.

The notion of pullback of a tensor field by a vector bundle map needs some attention.  Given a vector bundle map $\map{\Phi}{E}{E'}$ over a map $\map{\phi}{M}{M'}$, for every section $\beta$ of the $p$-covariant tensor bundle $\otimes ^pE'{}^*\to M'$ we define the section $\Phi\pb\beta$ of $\otimes^p E^*\to M$ by
$$
(\Phi\pb\beta)_m(a_1,a_2,\ldots,a_p)=
\beta_{\phi(m)}\bigl(\Phi(a_1),\Phi(a_2),\ldots,\Phi(a_p)\bigr).
$$
The tensor $\Phi\pb\beta$ is said to be the pullback of $\beta$ by $\Phi$. For a function $f\in\cinfty{M'}$ (i.e. for $p=0$) we just set $\Phi\pb f=\phi^*f=f\circ\phi$.  It follows that $\Phi\pb(\alpha\otimes\beta)=(\Phi\pb\alpha)\otimes(\Phi\pb\beta)$. 

When $\Phi$ is fiberwise invertible, we define the pullback of a section $\sigma$ of $E'$ by $(\Phi\pb\sigma)_m=(\Phi_m)^{-1}(\sigma(\phi(m))$, where $\Phi_m$ is the restriction of $\Phi$ to the fiber $E_m$ over $m\in M$. Thus, the pullback of any tensor field is defined.

In the case of the tangent bundles $E=TM$ and $E'=TM'$, and the tangent map $\map{T\varphi}{TM}{TM'}$ of a map $\map{\varphi}{M}{M'}$ we have that $(T\varphi)\pb\beta=\varphi^*\beta$ is the ordinary pullback by $\varphi$ of a tensorfield $\beta$ on $M'$. Notice the difference between the symbols ${}\pb$ (star) and ${}^*$ (asterisque). 

%=====================================================================

\section{Preliminaries on Lie algebroids}
\label{preliminaries}

In this section we recall some well known notions and introduce a few
concepts concerning the geometry of Lie algebroids. We refer the reader to~\cite{CaWe,Mackenzie} for
details about Lie groupoids, Lie algebroids and their role in
differential geometry.

\subsection*{Lie algebroids}
Let $M$ be an $n$-dimensional manifold and let $\map{\tau}{E}{M}$ be a vector bundle. A vector bundle map $\map{\rho}{E}{TM}$ over the identity is called an anchor map. The vector bundle $E$ together with an anchor map $\rho$ is said to be an \emph{anchored vector bundle}. A structure of \emph{Lie algebroid} on $E$ is given by a Lie algebra structure on the $\cinfty{M}$-module of sections of the bundle, $(\sec{E},[\cdot\,,\cdot])$, together with an anchor map, satisfying the compatibility condition
$$
[\sigma,f\eta] = f[\sigma,\eta] + \bigl(
\rho(\sigma)f \bigr) \eta .
$$
Here $f$ is a smooth function on $M$, $\sigma$, $\eta$ are
sections of $E$ and we have denoted by $\rho(\sigma)$ the vector field
on $M$ given by $\rho(\sigma)(m)=\rho(\sigma(m))$. From the compatibility
condition and the Jacobi identity, it follows that the map
$\sigma\mapsto\rho(\sigma)$ is a Lie algebra homomorphism from
$\sec{E}$ to $\mathfrak{X}(M)$. 

It is convenient to think of a Lie algebroid as a substitute of the tangent bundle of $M$. In this way, one regards an element $a$ of $E$ as a generalized velocity, and the actual velocity $v$ is obtained when applying the anchor to $a$, i.e., $v=\rho(a)$.
A curve $\map{a}{[t_0,t_1]}{E}$ is said to be \emph{admissible} if $\dot{m}(t)=\rho(a(t))$, where $m(t)=\tau(a(t))$ is the base curve.

When $E$ carries a Lie algebroid structure, the image of the anchor map, $\rho(E)$, defines an integrable smooth generalized distribution on $M$.  Therefore, $M$ is foliated by the integral leaves of $\rho(E)$, which are called the leaves of the Lie algebroid. It follows that $a(t)$ is admissible if and only if the curve $m(t)$ lies on a leaf of the Lie algebroid, and that two points are in the same leaf if and only if they are connected by (the base curve of) an admissible curve.

%No hace falta ??
A Lie algebroid is said to be transitive if it has only one leaf, which is obviously equal to $M$. It is easy to see that $E$ is transitive if and only if $\rho$ is surjective. If $E$ is not transitive, then the restriction of the Lie algebroid to a leaf $L\subset M$, $E_{|L}\to L$ is transitive.

Given local coordinates $(x^i)$ in the base manifold $M$ and a local basis $\{e_\alpha\}$ of sections of $E$, we have local coordinates $(x^i,y^\alpha)$ in $E$; if $a\in E_m$ is an element then we can write $a=a^\alpha e_\alpha(m)$ and thus the coordinates of $a$ are $(m^i,a^\alpha)$, where $m^i$ are the coordinates of the point $m$. The anchor map is locally determined by the local functions $\rho^i_\alpha$ on $M$ defined by $\rho(e_\alpha)=\rho^i_\alpha(\partial/\partial x^i)$. In addition, for a Lie algebroid, the Lie bracket is determined by the functions  $C^\alpha_{\beta\gamma}$ defined by $[e_\alpha,e_\beta]=C^\gamma_{\alpha\beta}e_\gamma$. The functions $\rho^i_\alpha$ and $C^\alpha_{\beta\gamma}$ are said to be the structure functions of the Lie algebroid in this coordinate system. They satisfy the following relations
\begin{align*}
  \rho^j_\alpha\pd{\rho^i_\beta}{x^j} -
  \rho^j_\beta\pd{\rho^i_\alpha}{x^j} = \rho^i_\gamma
  C^\gamma_{\alpha\beta},
  \quand
  \sum_{\mathrm{cyclic}(\alpha,\beta,\gamma)} \left[\rho^i_\alpha\pd{
      C^\nu_{\beta\gamma}}{x^i} + C^\mu_{\alpha\nu}
    C^\nu_{\beta\gamma}\right]=0.
\end{align*}
which are called the structure equations of the Lie algebroid.

\subsection*{Exterior differential}
The anchor $\rho$ allows us to define the differential of a function on the base manifold with respect to an element $a\in E$. It is given by
$$
df(a)=\rho(a)f.
$$
It follows that the differential of $f$ at the point $m\in M$ is an element of $E_m^*$. 

Moreover, a structure of Lie algebroid on $E$ allows us to extend the differential to sections of the bundle $\ext[p]{E}$ which we will call just $p$-forms. If $\omega\in\sec{\ext[p]{E}}$ then $d\omega\in\sec{\ext[p+1]{E}}$ is defined by 
\begin{align*}
d\omega(\sigma_0,\sigma_1,\ldots,\sigma_p)
&=
\sum_i\rho(\sigma_i)(-1)^i
  \omega(\sigma_0,\ldots,\widehat{\sigma_i},\ldots,\sigma_p)+\\
&\qquad{}+
\sum_{i<j}(-1)^{i+j}
\omega([\sigma_i,\sigma_j],\sigma_0,\ldots,\widehat{\sigma_i},\ldots,\widehat{\sigma_j},\ldots,\sigma_p).
\end{align*}
It follows that $d$ is a cohomology operator, that is $d^2=0$.

Locally the exterior differential is determined by
$$
dx^i=\rho^i_\alpha e^\alpha
\qquand
de^\gamma=-\frac{1}{2}C^\gamma_{\alpha\beta}e^\alpha\wedge e^\beta.
$$
In this paper the symbol $d$ will always denote the exterior differential on the Lie algebroid $E$ and not the ordinary exterior differential on a manifold $M$ to which it reduces when $E$ is the standard Lie algebroid $E=TM$ over $M$.

The usual Cartan calculus extends to the case of Lie algebroids (See~\cite{Nijenhuis}). For every section $\sigma$ of $E$ we have a derivation $i_\sigma$ (contraction) of degree $-1$ and a derivation $d_\sigma=i_\sigma\circ d+d\circ i_\sigma$ (Lie derivative) of degree $0$. Since $d^2=0$ we have that $d_\sigma\circ d=d\circ d_\sigma$.

\subsection*{Admissible maps and morphisms}
Let $\map{\tau}{E}{M}$ and $\map{\tau'}{E'}{M'}$ be two
anchored vector bundles, with anchor maps $\map{\rho}{E}{TM}$ and
$\map{\rho'}{E'}{TM'}$.  Let $\map{\Phi}{E}{E'}$ be a fiberwise linear map over $\map{\phi}{M}{M'}$. We will say that $\Phi$ is admissible if it maps admissible curves into admissible curves. It follows that $\Phi$ is admissible if and only if 
$T\phi\circ\rho = \rho'\circ\Phi$. This condition can be conveniently expressed in terms of the exterior differential as follows. The map $\Phi$ is \emph{admissible} if and only if $\Phi\pb df=d\Phi\pb f$ for every function $f\in\cinfty{M}$. If $E$ and $E'$ are Lie algebroids, then we say that $\Phi$ is a \emph{morphism} if $\Phi\pb d\theta=d\Phi\pb\theta$ for every $\theta\in\sec{\ext{E}}$. Obviously, a morphism is an admissible map. 

Let $(x^i)$ be a local coordinate system on $M$ and $(x'{}^i)$ a local coordinate system on $M'$. Let $\{e_\alpha\}$ and $\{e'_\alpha\}$ be local basis of sections of $E$ and $E'$, respectively, and $\{e^\alpha\}$ and $\{e'{}^\alpha\}$ the dual basis. The bundle map $\Phi$ is determined by the relations $\Phi\pb x'{}^i = \phi^i(x)$ and $\Phi\pb e'{}^\alpha = \phi^\alpha_\beta e^\beta$ for certain local functions $\phi^i$ and $\phi^\alpha_\beta$ on $M$. Then $\Phi$ is admissible if and only if
$$
  \rho^j_\alpha\pd{\phi^i}{x^j}=\rho'{}^i_\beta\phi^\beta_\alpha.
$$
The map $\Phi$ is a morphism of Lie algebroids if and only if 
$$
  \phi^\beta_\gamma C^\gamma_{\alpha\delta} =
    \left(\rho^i_\alpha\pd{\phi^\beta_\delta}{x^i} -
    \rho^i_\delta\pd{\phi^\beta_\alpha}{x^i}\right) +
  C'{}^\beta_{\theta\sigma}\phi^\theta_\alpha\phi^\sigma_\delta,
$$
in addition to the admissibility condition above. In these expressions $\rho^i_\alpha$, $C^\alpha_{\beta\gamma}$ are the structure functions on $E$ and $\rho'{}^i_\alpha$, $C'{}^\alpha_{\beta\gamma}$ are the structure functions on $E'$. 

\subsection*{Flow defined by a section} 
Every section of a Lie algebroid has an associated local flow composed of morphisms of Lie algebroids. 

\begin{theorem}
Let $\sigma$ be a section of a Lie algebroid $\map{\tau}{E}{M}$. There exists a local flow $\map{\Phi_s}{E}{E}$ such that 
$$
d_\sigma\theta=\frac{d}{ds}\Phi_s\pb\theta\at{s=0},
$$
for any section $\theta$ of $\ext{E}$. For every fixed $s$, the map $\Phi_s$ is a morphism of Lie algebroids. The base maps $\map{\phi_s}{M}{M}$ are the flow of the vectorfield $\rho(\sigma)\in\vectorfields{M}$. 
\end{theorem}
\begin{proof}
Indeed, let $X_\sigma\spC$ be the vector field on $E$ determined by its action on fiberwise linear functions
$X_\sigma\spC\hat{\theta}=\widehat{d_\sigma\theta}$, for every section $\theta$ of $E^*$, and which is known as the complete lift of $\sigma$ (see~\cite{GrUr,LMLA}). It follows that $X_\sigma\spC\in\vectorfields{E}$ is projectable and projects to the vectorfield $\rho(\sigma)\in\vectorfields{M}$. Thus, if we consider the flow $\Phi_s$ of $X_\sigma\spC$ and the flow $\phi_s$ of $\rho(\sigma)$, we have that $\Phi_s$ is a bundle map over $\phi_s$. Moreover, since $X_\sigma\spC$ is linear (maps linear functions into linear functions) we have that $\Phi_s$ is a vector bundle map for every $s$. 

In order to prove the relation $d_\sigma\theta=\frac{d}{ds}\Phi_s\pb\theta\at{s=0}$, it is sufficient to prove it for 1-forms $\theta$. Then for every $m\in M$ and $a\in E_m$ we have 
\begin{align*}
\pai{(d_\sigma\theta)_m}{a}
&=\widehat{d_\sigma\theta}(a)
=(X_\sigma^C\hat{\theta})(a)
=\frac{d}{ds}(\hat{\theta}\circ\Phi_s)\at{s=0}(a)\\
&=\frac{d}{ds}\pai{\theta_{\phi_s(m)}}{\Phi_s(a)}\at{s=0}
=\frac{d}{ds}\Phi_s\pb\theta\at{s=0}(a).
\end{align*}

Moreover, the maps $\Phi_s$ are morphisms of Lie algebroids. Indeed, from the relation $d\circ d_\sigma=d_\sigma\circ d$ we have that $\frac{d}{ds}(\Phi_s\pb d\theta-d\Phi_s\pb\theta)\at{s=0}=0$. It follows that $\Phi_s\pb d\theta=d\Phi_s\pb\theta$, so that $\Phi_s$ is a morphism.
\end{proof}

By duality, it follows that 
$
d_\sigma\zeta=\frac{d}{ds}\Phi_s\pb\zeta\big|_{s=0}
$ 
for every section $\zeta$ of $E$.  Therefore we have a similar formula for the derivative of any tensorfield.

\begin{remark}
In the case of the standard Lie algebroid $E=TM$ we have that a section $\sigma$ is but a vector field on $M$. In this case, the vector field $X_\sigma\spC$ is but the usual complete lift of the vectorfield $\sigma$ and the flow defined by the section $\sigma$ is but $\Phi_s=T\phi_s$ where $\phi_s$ is the  flow of the vectorfield $\sigma$.
\end{remark}

\section{Jets}
\label{jetoids}

We consider two vector bundles $\map{\prEM}{E}{M}$ and $\map{\prFN}{F}{N}$ and a surjective vector bundle map $\map{\prEF}{E}{F}$ over the map $\map{\prMN}{M}{N}$.
Moreover, we will assume that $\map{\prMN}{M}{N}$ is a smooth fiber bundle. We will denote by $K\to M$ the kernel of the map $\prEF$, which is a vector bundle over~$M$. Given a point $m\in M$, if we denote $n=\prMN(m)$, we have the following exact sequence
$$
 0\to K_m \to E_m \to F_n \to 0,
$$
and we can consider the set $\Jpi[m]$ of splittings $\phi$ of such sequence. More concretely, we define the following sets 
\begin{align*}
\Lpi[m]&=\set{\map{w}{F_n}{E_m}}{\text{$w$ is linear}}\\
\Jpi[m]&=\set{\phi\in\Lpi[m]}{\prEF\circ\phi=\id_{F_n}}\\
\Vpi[m]&=\set{\psi\in\Lpi[m]}{\prEF\circ\psi=0}.
\end{align*}
Therefore $\Lpi[m]$ is a vector space, $\Vpi[m]$ is a vector subspace of $\Lpi[m]$ and $\Jpi[m]$ is an affine subspace of $\Lpi[m]$ modeled on the vector space is $\Vpi[m]$. By taking the union, $\Lpi=\cup_{m\in M}\Lpi[m]$, $\Jpi=\cup_{m\in M}\Jpi[m]$ and $\Vpi=\cup_{m\in M}\Vpi[m]$, we get the vector bundle $\map{\prLM}{\Lpi}{M}$ and the affine subbundle $\map{\prJM}{\Jpi}{M}$ modeled on the vector bundle $\map{\prVM}{\Vpi}{M}$. We will also consider the projection $\map{\pi_1}{\Jpi}{N}$ defined by composition $\pi_1=\pi\circ\pi_{10}$.

\begin{remark}
In the standard case, one has a bundle $\map{\prMN}{M}{N}$ and then considers the tangent spaces at $m\in M$ and $n=\prMN(m)\in N$ together with the differential of the projection $\map{T_m\prMN}{T_mM}{T_nN}$. A 1-jet at $m$, in the standard sense, is a linear map $\map{\phi}{T_nN}{T_mM}$ such that $T_m\prMN\circ\phi=\id_{T_nN}$. It follows that there exist sections $\varphi$ of the bundle $\prMN$ such that $\varphi(n)=m$ and $T_n\varphi=\phi$. Thus, in the case of the tangent bundles $\map{\tau_N}{F\equiv TN}{N}$ and $\map{\tau_M}{E\equiv TM}{M}$ with $\map{\prEF=T\prMN}{TM}{TN}$, our definition is equivalent to the standard definition of a 1-jet at a point $m$ as an equivalence class of sections which has the same value and same first derivative at the point $m$. With the standard notations~\cite{Saunders}, we have that $J^1\nu\equiv\mathcal{J}(T\nu)$. Obviously, this example leads our developments throughout this paper. 
\end{remark}

In view of this fact, an element of $\Jpi[m]$ will be simply called a \emph{jet} at the point $m\in M$ and accordingly the bundle $\Jpi$ is said to be the first \emph{jet bundle} of $\prEF$.

\medskip

Local coordinates on $\Jpi$ are given as follows. We consider local coordinates $(\bar{x}^i)$ on $N$ and $(x^i,u^A)$ on $M$ adapted to the projection $\prMN$, that is $\bar{x}^i\circ\prMN=x^i$. We also consider local basis of sections $\{\bar{e}_a\}$ of $F$ and $\{e_a,e_\alpha\}$ of $E$ adapted to the projection $\prEF$, that is $\prEF\circ e_a=\bar{e}_a\circ\prMN$ and  $\prEF\circ e_\alpha=0$. In this way $\{e_\alpha\}$ is a base of sections of $K$. An element $w$ in $\Lpi[m]$ is of the form $w=(w_a^be_b+w_a^\alpha e_\alpha)\otimes \bar{e}^a$, and it is in $\Jpi[m]$ if and only if $w_a^b=\delta_a^b$, i.e.\ an element $\phi$ in $\Jpi$ is of the form $\phi=(e_a+\phi^\alpha_a e_\alpha)\otimes\bar{e}_a$. If we set $y^\alpha_a(\phi)=\phi^\alpha_a$, we have adapted local coordinates $(x^i,u^A,y^\alpha_a)$ on $\Jpi$. 
Similarly, an element $\psi\in\Vpi[m]$ is of the form $\psi=\psi^\alpha_a e_\alpha\otimes e^a$. If we set $y^\alpha_a(\psi)=\psi^\alpha_a$, we have adapted local coordinates $(x^i,u^A,y^\alpha_a)$ on $\Vpi$. As usual, we use the same name for the coordinates in an affine bundle and in the associated vector bundle. 

An element $\zeta\in\dLpi[m]$ defines an affine function $\hat{z}$ on $\Jpi[m]$ by contraction $\hat{z}(\phi)=\pai{z}{\phi}$ where $\pai{\cdot}{\cdot}$ is the pairing $\pai{z}{\phi}=\tr(z\circ\phi)=\tr(\phi\circ z)$. Therefore, a section $\theta$ of $\dLpi$ defines a fiberwise affine function $\hat{\theta}$ on $\Jpi$, 
$
\hat{\theta}(\phi)=\pai{\theta_{\prJM(\phi)}}{\phi}=
\tr(\theta_{\prJM(\phi)}\circ\phi).
$
In local coordinates, a section of $\dLpi$ is of the form $\theta=(\theta^a_b(x)e^b+\theta^a_\alpha(x)e^\alpha)\otimes\bar{e}_a$, and the affine function defined by $\theta$ is
$$
\hat{\theta}=\theta^a_a(x)+\theta^a_\alpha(x)y^\alpha_a.
$$

\subsection*{Anchor}
Consider now anchored structures on the bundles $E$ and $F$, that is, we have two vector bundle maps $\map{\rho_F}{F}{TN}$ and $\map{\rho_E}{E}{TM}$ over the identity in $N$ and $M$ respectively. We will assume that the map $\pi$ is admissible, that is $\rho_F\circ\prEF=T\prMN\circ\rho_E$. Therefore we have
$$
\rho_F(\bar{e}_a)=\rho^i_a\pd{}{\bar{x}^i}
\qquand
\left\{\begin{array}{l}
\rho_E(e_a)=\rho^i_a\dpd{}{x^i}+\rho^A_a\dpd{}{u^A}\\[5pt]
\rho_E(e_\alpha)=\rho^A_\alpha\dpd{}{u^A},
\end{array}\right.
$$
with $\rho^i_a=\rho^i_a(x)$, $\rho^A_a=\rho^A_a(x,u)$ and $\rho^A_\alpha=\rho^A_\alpha(x,u)$.
Equivalently
$$
d\bar{x}^i=\rho^i_a\bar{e}^a
\qquand
\left\{\begin{array}{l}
dx^i=\rho^i_a e^a\\
du^A=\rho^A_a e^a+\rho^A_\alpha e^\alpha.
\end{array}\right.
$$

\medskip

The anchor allows us to define the concept of total derivative of a function with respect to a section. Given a section $\sigma\in\sec{F}$, the total derivative of a function $f\in\cinfty{M}$ with respect to $\sigma$ is the function $\widehat{df\otimes\sigma}$, i.e. the affine function associated to $df\otimes\sigma\in\sec{\dLpi}$. In particular, the total derivative with respect to an element $\bar{e}_a$ of the local basis of sections of $F$,  will be denoted by $\p{f}_{|a}$. In this way, if $\sigma=\sigma^a\bar{e}_a$ then 
$
\widehat{df\otimes\sigma}=\p{f}_{|a}\sigma^a,
$
where the coordinate expression of $\p{f}_{|a}$ is
$$
\p{f}_{|a}=\rho^i_a\pd{f}{x^i}+(\rho^A_a+\rho^A_\alpha y^\alpha_a)\pd{f}{u^A}.
$$
Notice that, for a function $f$ in the base $N$, we have that $\p{f}_{|a}=\rho^i_a\pd{f}{x^i}$ are just the components of $df$ in the basis $\{\bar{e}_a\}$.

\subsection*{Bracket}

Let us finally assume that we have Lie algebroid structures on $\map{\prFN}{F}{N}$ and on $\map{\prEM}{E}{M}$, and that the projection $\prEF$ is a morphism of Lie algebroids. This condition implies the vanishing of some structure functions. 

We have the following expressions for the  brackets of elements in the basis of sections 
$$
[\bar{e}_a,\bar{e}_b]=C^c_{ab}\bar{e}_c
\qquand
\left\{\begin{aligned}
&[e_a,e_b]=C^\gamma_{ab}e_\gamma+C^c_{ab}e_c\\
&[e_a,e_\beta]=C^\gamma_{a\beta}e_\gamma \\
&[e_\alpha,e_\beta]=C^\gamma_{\alpha\beta}e_\gamma.
 \end{aligned}\right.
$$
where $C^a_{bc}=C^a_{bc}(x)$ is a basic function. The exterior differentials in $F$ and $E$ are determined by
$$
d\bar{e}^a=-\frac{1}{2}C^a_{bc}\bar{e}^b\wedge\bar{e}^c
\quand
\left\{\begin{array}{l}
de^a=-\frac{1}{2}C^a_{bc}e^b\wedge e^c\\[5pt]
de^\alpha=-\frac{1}{2}C^\alpha_{bc}e^b\wedge e^c
          -           C^\alpha_{b\gamma}e^b\wedge e^\gamma
          -\frac{1}{2}C^\alpha_{\beta\gamma}e^\beta\wedge e^\gamma.
 \end{array}\right.
$$

\subsubsection*{Affine structure functions}
We define the affine functions $Z^\alpha_{a\gamma}$ and $Z^\alpha_{ac}$ by
$$
Z^\alpha_{a\gamma}=\widehat{(d_{e_\gamma}e^\alpha)\otimes\bar{e}_a}
\qquand
Z^\alpha_{ac}=\widehat{(d_{e_c}e^\alpha)\otimes\bar{e}_a}.
$$
Explicitly, we have
$$
Z^\alpha_{a\gamma}=C^\alpha_{a\gamma}+C^\alpha_{\beta\gamma}y^\beta_a
\qquand
Z^\alpha_{ac}=C^\alpha_{ac}+C^\alpha_{\beta c}y^\beta_a.
$$
Indeed,
\begin{align*}
(d_{e_\gamma}e^\alpha)\otimes\bar{e}_a
&=i_{e_\gamma}\left(-\frac{1}{2}C^\alpha_{bc}e^b\wedge e^c-C^\alpha_{b\theta}e^b\wedge e^\theta-\frac{1}{2}C^\alpha_{\beta\theta}e^\beta\wedge e^\theta\right)
\otimes\bar{e}_a\\
&=(C^\alpha_{b\gamma}e^b+C^\alpha_{\beta\gamma}e^\beta)\otimes\bar{e}_a,
\end{align*}
and thus
$
Z^\alpha_{a\gamma}
=\widehat{(C^\alpha_{b\gamma}e^b+C^\alpha_{\beta\gamma}e^\beta)
  \otimes\bar{e}_a}
=C^\alpha_{a\gamma}+C^\alpha_{\beta\gamma}y^\beta_a
$.
The second coordinate expression can be found in a similar way.

%
%Similarly
%\begin{align*}
%(d_{e_c}e^\alpha)\otimes\bar{e}_a
%&=i_{e_c}\left(-\frac{1}{2}C^\alpha_{bd}e^b\wedge e^d-C^\alpha_{b\theta}e^b\wedge e^\theta-\frac{1}{2}C^\alpha_{\beta\theta}e^\beta\wedge e^\theta\right)
%\otimes\bar{e}_a\\
%&=(C^\alpha_{bc}e^b+C^\alpha_{\beta c}e^\beta)\otimes\bar{e}_a,
%\end{align*}
%and thus
%$$
%Z^\alpha_{ac}
%=\widehat{(C^\alpha_{bc}e^b+C^\alpha_{\beta c}e^\beta)
%  \otimes\bar{e}_a}
%=C^\alpha_{ac}+C^\alpha_{\beta c}y^\beta_a.
%$$

\section{Morphisms and admissible maps}
\label{morphisms}

By a section of $\pi$ we mean a vector bundle map $\Phi$ such that $\pi\circ\Phi=\id$, (i.e. we consider only linear sections of $\prEF$). It follows that the base map $\map{\phi}{N}{M}$ is a section of $\prMN$, i.e. $\prMN\circ\phi=\id_N$. The set of sections of $\pi$ will be denoted by $\sec{\pi}$. We will find local conditions for this map to be an admissible map between anchored vector bundles and also local conditions for this map to be a morphism of Lie algebroids. The set of those sections of $\pi$ which are a morphism of Lie algebroids will be denoted by $\Mor{\pi}$.

Taking adapted local coordinates $(x^i,u^A)$ on $M$, the map $\phi$ has the expression $\phi(x^i)=(x^i,\phi^A(x))$. If we moreover take an adapted basis $\{e_a,e_\alpha\}$ of local sections of $E$, then the expression of $\Phi$ is given by $\Phi(\bar{e}_a)=e^a+\phi^\alpha_ae_\alpha$, so that the map $\Phi$ is determined by the functions $\bigl(\phi^A(x),\phi^\alpha_a(x)\bigr)$ locally defined on $N$. The action on the dual basis is $\Phi\pb e^a=\bar{e}^a$, and $\Phi\pb e^\alpha=\phi^\alpha_a \bar{e}^a$, and for the coordinate functions $\Phi\pb x^i=\bar{x}^i$ and $\Phi\pb u^A=\phi^A$.

Let us see what the admissibility condition $\Phi\pb(df)=d(\Phi\pb f)$ means for this maps. We impose this condition to the coordinate functions. Taking $f=x^i$ we get an identity (i.e. no new condition arises), and taking $f=u^A$ we get
$$
d\phi^A
=d(\Phi\pb u^A)
=\Phi\pb(du^A)
=\Phi\pb(\rho^A_ae^a+\rho^A_\alpha e^\alpha)
=[(\rho^A_a\circ\phi)+(\rho^A_\alpha\circ\phi)\phi^\alpha_a] \, \bar{e}^\alpha,
$$
from where we get that $\Phi$ is an admissible map if and only if
$$
\rho^i_a\pd{\phi^A}{x^i}
=(\rho^A_a\circ\phi)+(\rho^A_\alpha\circ\phi)\phi^\alpha_a.
$$
In order to simplify the writing, we will apply the usual convention and we will omit the composition with $\phi$, since we know that this is an equation to be satisfied at the point $m=\phi(n)=(x^i,\phi^A(x))$ for every $n\in N$. With this convention, the above equation is written as
$$
\rho^i_a\pd{u^A}{x^i}=\rho^A_a+\rho^A_\alpha y^\alpha_a.
$$

Let us now see what the condition of being a morphism means in coordinates. If we impose $\Phi\pb de^a=d\Phi\pb e^a$ we get an identity, so that we just have to impose $\Phi\pb de^\alpha=d\Phi\pb e^\alpha$. On one hand we have
$$
d(\Phi\pb e^\alpha)
=d(\phi^\alpha_a\bar{e}^a)
=\frac{1}{2}\left(
  \rho^i_b\pd{\phi^\alpha_c}{x^i}-\rho^i_c\pd{\phi^\alpha_b}{x^i}
 -\phi^\alpha_a C^a_{bc}
  \right)\bar{e}^b\wedge\bar{e}^c
$$
and on the other hand
\begin{align*}
\Phi\pb d(e^\alpha)
&=-\Phi\pb\left(
   \frac{1}{2}C^\alpha_{\beta\gamma}e^\beta\wedge e^\gamma
  -C^\alpha_{b\gamma}e^b\wedge e^\gamma
  -\frac{1}{2}C^\alpha_{bc}e^b\wedge e^c
     \right)\\
&=-\frac{1}{2}\left(
   C^\alpha_{\beta\gamma}\phi^\beta_b\phi^\gamma_c
  +C^\alpha_{b\gamma}\phi^\gamma_c-C^\alpha_{c\gamma}\phi^\gamma_b
  -C^\alpha_{bc}\right)\bar{e}^b\wedge\bar{e}^c
\end{align*}
Thus, the bundle map $\Phi$ is a morphism if and only if it satisfies
$$
\rho^i_b\pd{\phi^\alpha_c}{x^i}-\rho^i_c\pd{\phi^\alpha_b}{x^i}
 -\phi^\alpha_a C^a_{bc}
+
   C^\alpha_{\beta\gamma}\phi^\beta_b\phi^\gamma_c
  +C^\alpha_{b\gamma}\phi^\gamma_c-C^\alpha_{c\gamma}\phi^\gamma_b
=C^\alpha_{bc},
$$
in addition to the admissibility condition. Taking into account our notation $\p{f}_{|a}=\rho^i_a\pd{f}{x^i}$ for a function $f\in\cinfty{N}$, we can write the above expressions as
\begin{gather*}
\p{u}^A_{|a}=\rho^A_a+\rho^A_\alpha y^\alpha_a\\
%\intertext{and}
 \bigr(\p{y}^\alpha_{c|b}+C^\alpha_{b\gamma}y^\gamma_c\bigr)
-\bigr(\p{y}^\alpha_{b|c}+C^\alpha_{c\gamma}y^\gamma_b\bigr)
 +C^\alpha_{\beta\gamma}y^\beta_by^\gamma_c-y^\alpha_a C^a_{bc}
=C^\alpha_{bc}.
\end{gather*}
where we recall that they are to be satisfied at every point $m=\phi(n)$.

\section{Complete lift of a section}
\label{complete.lift}

In this section we will define the lift of a projectable section of $E$ to a vectorfield on $\Jpi$, in a similar way to the definition of the  first jet prolongations of a projectable vectorfield in the standard theory of jet bundles~\cite{Saunders}.

We consider a section $\sigma$ of a Lie algebroid $E$ projectable over a section $\eta$ of $F$. We denote by $\Psi_s$ the flow on $E$ associated to $\sigma$ and by $\Phi_s$ the flow on $F$ associated to $\eta$. We recall that, for every fixed $s$, the maps $\Psi_s$ and $\Phi_s$ are morphisms of Lie algebroids. Moreover, the base maps $\psi_s$ and $\phi_s$, are but the flows of the vectorfields $\rho_E(\sigma)$ and $\rho_F(\eta)$, respectively.

The projectability of the section implies the projectability of the flow. It follows that (locally, in the domain of the flows) we have defined a map $\map{\Lprol{\Psi_s}}{\Lpi}{\Lpi}$ by means of 
$$
\Lprol{\Psi_s}(w)=\Psi_s\circ w\circ\Phi_{-s}.
$$
By restriction of $\Lprol{\Psi_s}$ to $\Jpi$ we get a map $\Jprol{\Psi_s}$, which is a local flow in $\Jpi$. We will denote by $\XsigmaJ$ the vector field on $\Jpi$ generating the flow $\Jprol{\Psi_s}$. The vectorfield $\XsigmaJ$ will be called the \emph{complete lift} to $\Jpi$ of the section $\sigma$. Since $\Jprol{\Psi_s}$ projects to the flow $\Psi_s$ it follows that the vector field $\XsigmaJ$ projects to the vector field $\rho^E(\sigma)$ in $M$. 

\subsection*{Derivative of a section of $\dLpi$}
In order to find a more operational redefinition of the complete lift, we consider the derivative of a section of $\dLpi$ with respect to a projectable section. If we are given a projectable section $\sigma$ of $E$, then we can define the Lie derivative of $\theta\in\sec{\dLpi}$ with respect to $\sigma$ by means of 
$$
d_\sigma\theta
=\frac{d}{ds}\Psi_{s}\pb\theta\at{s=0}.
$$
Here by $\Psi_s\pb\theta$ we mean 
$(\Psi_s\pb\theta)_m=\Phi_{-s}\circ\theta_{\psi_s(m)}\circ\Psi_s$, so that the above derivative is given explicitly by
$$
(d_\sigma\theta)(a)
=\frac{d}{ds}
\Bigl(\Phi_{-s}\bigl(\theta_{\psi_s(m)}(\Psi_s(a))\bigr)\Bigr)
\at{s=0},
$$
for every $a\in E_m$. It can be easily seen that this prescription defines a derivation of the module of sections of $\dLpi$, with associated vector field $\rho^E(\sigma)$. 

An alternative definition can be given as follows. Given $\theta\in\sec{\dLpi}$ we take $P$ any section of $E^*\otimes E$ such that $\theta=\prEF\circ P$. Then we have that
$$
d_\sigma\theta=\prEF\circ d_\sigma P.
$$
The right hand side is independent of the choice of $P$ because $K$ is a $\prEF$-ideal. Indeed, if we take two different tensors $P$ they differ by a section $\zeta$ of $E^*\otimes K$. Then, for any section $\eta\in\sec{E}$,
$$
\prEF\bigl((d_\sigma\zeta)(\eta)\bigr)
=\prEF\bigl(d_\sigma(\zeta(\eta))\bigr)- \prEF\bigl(\zeta(d_\sigma\eta)\bigr)
=\prEF\bigl(d_\sigma(\zeta(\eta))\bigr),
$$
which vanishes, because $\sigma$ is projectable and $\theta(\eta)$ a section of $K$.

The equivalence of both definitions follows by taking into account that 
$(\Psi_s\pb P)_m=\Psi_{-s}\circ P_{\psi(m)}\circ\Psi_s$ and thus
\begin{align*}
(\Psi_s\pb\theta)_m
&=\Phi_{-s}\circ\theta_{\psi(m)}\circ\Psi_s
=\Phi_{-s}\circ\prEF\circ P_{\psi(m)}\circ\Psi_s\\
&=\prEF\circ\Psi_{-s}\circ P_{\psi(m)}\circ\Psi_s
=\prEF\circ(\Psi_s\pb P)_m.
\end{align*}
In particular, the above definition implies that the usual rules for calculating Lie derivatives apply to the calculation of $d_\sigma\theta$. 

Notice the following relation between the action of $\Lprol{^*\Psi_s}$ and the pullback 
$$
(\Psi_s\pb\theta)_m=\mathcal{L}^*\Psi_s(\theta_{\psi_s(m)}).
$$
Indeed, for every $m\in M$ and $w\in\Lpi[m]$ we have
\begin{align*}
\pai{\Lprol{^*\Psi_s}(\theta_{\psi_s(m)})}{w}
&=\pai{\theta_{\psi_s(m)}}{\Lprol{\Psi_s}(w)}
=\tr(\theta_{\psi_s(m)}\circ\Psi_s\circ w\circ\Phi_{-s})\\
&=\tr(\Phi_{-s}\circ\theta_{\psi_s(m)}\circ\Psi_s\circ w)
=\pai{(\Psi_s\pb\theta)_m}{w}.
\end{align*}

With the help of the derivative of a section of $\dLpi$ we can characterize the complete lift of a projectable section of $E$ in terms of its action on affine functions as follows.

\begin{proposition}
Given a projectable section $\sigma$ of $E$ the vectorfield $\XsigmaJ$ is is characterized by the following properties
\begin{itemize}
\item $\XsigmaJ$ is projectable and projects to $\rho(\sigma)$, and
\item $\mathcal{L}_\XsigmaJ\hat{\theta}=\widehat{d_\sigma\theta}$, for every section $\theta$ of $\dLpi$.
\end{itemize}
\end{proposition}
\begin{proof}
It is easy to see that both conditions are compatible, and therefore they define a unique vectorfield in $\Jpi$. We just have to prove that $\XsigmaJ$ satisfies those properties. From the definition of $\XsigmaJ$ it is clear that it projects to $\rho(\sigma)$. Moreover, for every section $\theta$ of $\Lpi^*$ we have that
\begin{align*}
\hat{\theta}(\Jprol{\Psi_s}(\phi))
&=\pai{\theta_{\psi_s(m)}}{\Jprol{\Psi_s}(\phi)}
=\pai{\theta_{\psi_s(m)}}{\Lprol{\Psi_s}(\phi)}\\
&=\pai{\Lprol{^*\Psi_s}\bigl(\theta_{\psi_s(m)}\bigr)}{\phi}
=\pai{(\Psi_s\pb\theta)_m}{\phi}
\end{align*}
and taking the derivative with respect to $s$ at $s=0$,
$$
X^{\Jpi}_\sigma(\hat{\theta})(\phi)
=\frac{d}{ds}\hat{\theta}(\Jprol{\Psi_s}(\phi))\Big|_{s=0}
=\frac{d}{ds}\pai{(\Psi_s\pb\theta)_m}{\phi}\Big|_{s=0}
=\pai{(d_\sigma\theta)_m}{\phi}\\
=\widehat{d_\sigma\theta}(\phi)
$$
which concludes the proof.
\end{proof}

\subsubsection*{Local expression}
Locally, a section $\sigma=\sigma^a e_a+\sigma^\alpha e_\alpha$ is projectable if $\sigma^a=\sigma^a(x^i)$ depends only on $x^i$. If its complete lift $\XsigmaJ$ projects to $\rho(\sigma)$ it must be of the form
$$
\XsigmaJ=\rho^i_a\sigma^a\pd{}{x^i}+
(\rho^A_a\sigma^a+\rho^A_\alpha\sigma^\alpha)\pd{}{u^A}+\sigma^\alpha_a\V^a_\alpha.
$$
where the function $\sigma^\alpha_a$ are to be determined by the second condition. If we take the local section $\theta=e^\alpha\otimes\bar{e}_a$ of $\dLpi$, then $\hat{\theta}=y^\alpha_a$, and thus $\sigma^\alpha_a
=d_{\sigma\spC}y^\alpha_a
=d_{\sigma\spC}\hat{\theta}
=\widehat{d_{\sigma}\theta}
$. We have
$$
d_\sigma\theta =
(d_\sigma e^\alpha)\otimes\bar{e}_a+e^\alpha\otimes(d_\eta\bar{e}_a)=
(d_\sigma e^\alpha)\otimes\bar{e}_a-e^\alpha\otimes d_{\bar{e}_a}\eta.
$$
Also
$$
d_\sigma e^\alpha=d i_\sigma e^\alpha+i_\sigma de^\alpha
=d\sigma^\alpha+\sigma^bi_{e_a} de^\alpha+\sigma^\beta i_{e_\beta} de^\alpha.
$$
and
$$
d_{\bar{e}_a}\eta=[\bar{e}_a,\eta]=
\bigl(\rho^i_a\pd{\sigma^b}{x^i}+\sigma^cC^b_{ac}\bigr)\bar{e}_b
$$
Putting all together
$$
d_\sigma\theta
= d\sigma^\alpha\otimes\bar{e}_a
 + \sigma^b(i_{e_b} de^\alpha\otimes\bar{e}_a)
 + \sigma^\beta(i_{e_\beta} de^\alpha\otimes\bar{e}_a)
 -\Bigl(\p{\sigma}^b_{|a}+\sigma^cC^b_{ac}\Bigr)
    e^\alpha\otimes
 \bar{e}_b,
$$
and taking into account the definition of the functions $Z^\alpha_{ab}$ and $Z^\alpha_{a\beta}$ we arrive to
$$
\sigma^\alpha_a
=
\p{\sigma}^\alpha_{|a}+Z^\alpha_{ab}\sigma^b+Z^\alpha_{a\beta}\sigma^\beta
-y^\alpha_b\Bigl(\p{\sigma}^b_{|a}+\sigma^cC^b_{ac}\Bigr).
$$
In particular, if $\sigma$ projects to $0$, i.e.\ $\sigma^a=0$, we have
$$
\XsigmaJ=\rho^A_\alpha\sigma^\alpha\pd{}{u^A}+\left(
  \p{\sigma}^\alpha_{|a}+Z^\alpha_{a\beta}\sigma^\beta
  \right)\pd{}{y^\alpha_a}.
$$

\section{Variational Calculus}
\label{variational}

In what follows in this paper we consider the case where the Lie algebroid $F$ is the tangent bundle $F=TN$ with $\rho_F=\id_{TN}$ and $[\cdot,\cdot]$ the usual Lie bracket of vectorfields on $N$. The Lie algebroid $E$ remains a general Lie algebroid.

\subsection*{Variational problem}
Given a Lagrangian function $L\in\cinfty{\Jpi}$ and a volume form $\omega\in\ext[r]{(TN)}$, where $r=\operatorname{dim}(N)$, we consider the following variational problem: find the critical points of the action functional 
$$
\CMcal{S}(\Phi)=\int_N L(\Phi)\,\omega
$$
defined on the set of sections of $\pi$ which are moreover morphisms of Lie algebroids, that is, defined on the set $\Mor{\prEF}$. Here by $L(\Phi)$ we mean the function $n\mapsto L(\Phi_n)$, where $\Phi_n\in\Jpi$ is the restriction of $\Phi$ the fiber $F_n=T_nN$.

It is important to notice that the above variational problem is a constrained problem, not only because the condition $\pi\circ\Phi=\id_F$, which can be easily solved, but because of the condition of $\Phi$ being a morphism of Lie algebroids, which is a condition on the derivatives of $\Phi$. Taking coordinates on $N$ such that the volume form is $\omega=dx^1\wedge \cdots\wedge dx^r$, the problem is to
find the critical points of 
\begin{align*}
&\int_N L(x^i,u^A,y^\alpha_a)\,dx^1\wedge \cdots\wedge dx^r\\
\intertext{subject to the constraints}
&\pd{u^A}{x^a}=\rho^A_a+\rho^A_\alpha y^\alpha_a\\
&\pd{y^\alpha_c}{x^b}-\pd{y^\alpha_b}{x^c}+C^\alpha_{b\gamma}y^\gamma_c
-C^\alpha_{c\gamma}y^\gamma_b
 +C^\alpha_{\beta\gamma}y^\beta_by^\gamma_c-y^\alpha_a C^a_{bc}
=C^\alpha_{bc}.
\end{align*}

The first method one can try for solving the problem is to use Lagrange multipliers. Nevertheless, one has no warranties that all solutions to this problem are normal (i.e. not strictly abnormal). In fact, in very simple cases as the problem of a Heavy top~\cite{LMLA} one can easily see that there will be strictly abnormal solutions. Therefore we will try another alternative, which consists in finding explicitly finite variations of a solution, that is, by defining a curve in $\Mor{\pi}$ starting at the given solution.

\subsection*{Equations for critical sections}
In order to find admissible variations we consider sections of $E$ and the associated flow. With the help of this flow we can transform morphisms of Lie algebroids into morphisms of Lie algebroids.

Let $\Phi$ be a critical point of $\CMcal{S}$. In order to find admissible variations we consider a $\pi$-vertical section $\sigma$ of $E$. Its flow $\map{\Psi_s}{E}{E}$ projects to the identity in $F=TN$. Moreover we will require $\sigma$ to have compact support. Since for every fixed $s$, the map $\Psi_s$ is a morphism of Lie algebroids, it follows that the map $\Phi_s=\Psi_s\circ\Phi$ is a section of $\pi$ and a morphism of Lie algebroids, that is, $s\mapsto\Phi_s$ is a curve in $\mathcal{M}(\pi)$. Using this kind of variations we have the following result.

\begin{theorem}
A map $\Phi$ is a critical section of $\CMcal{S}$ if and only if in local coordinates such that the volume form is $\omega=dx^1\wedge\cdots\wedge dx^r$ the components $y^\alpha_a$ of  $\Phi$ satisfy the system of partial differential equations
\begin{align*}
&\pd{u^A}{x^a}=\rho^A_a+\rho^A_\alpha y^\alpha_a\\
& \pd{y^\alpha_a}{x^b}-\pd{y^\alpha_b}{x^a}
 +C^\alpha_{b\gamma}y^\gamma_a-C^\alpha_{a\gamma}y^\gamma_b
 +C^\alpha_{\beta\gamma}y^\beta_by^\gamma_a
 +C^\alpha_{ab}=0\\
&\frac{d\,\,}{dx^a}\left(\pd{L}{y^\alpha_a}\right)
=\pd{L}{y^\gamma_a}Z^\gamma_{a\alpha}
+\pd{L}{u^A}\rho^A_\alpha.
\end{align*}
\end{theorem}
\begin{proof}
Recall that by $L(\Phi)$ we mean the function in $N$ given by $L(\Phi)(n)=L(\Phi_n)$, where $\Phi_n$ is the restriction of $\map{\Phi}{F}{E}$ to the fibre $F_n$. The function $L(\Phi_s)$ is 
$$
L(\Phi_s)(n)=L(\Psi_s\circ\Phi_n)=L(\Jprol{\Psi_s}(\Phi_n)
=(\Jprol{\Psi_s}^*L)(\Phi)(n),
$$
and therefore the variation of the action along the curve $s\mapsto\Phi_s$ is
$$
0=\frac{d}{ds}\CMcal{S}(\Phi_s)\at{s=0}
=\int_N\frac{d}{ds}L(\Phi_s)\at{s=0}\omega
=\int_N(\mathcal{L}_\XsigmaJ L)(\Phi)\,\omega.
$$
Taking into account the local expression of $\XsigmaJ$ for a $\pi$-vertical $\sigma$, we have that 
\begin{align*}
\mathcal{L}_\XsigmaJ L
&=\rho^A_\alpha\sigma^\alpha\pd{L}{u^A}+\left(
  \frac{d\sigma^\alpha}{dx^a}+Z^\alpha_{a\beta}\sigma^\beta
  \right)\pd{L}{y^\alpha_a}\\
&=\sigma^\alpha\left[\rho^A_\alpha\pd{L}{u^A}+ Z^\gamma_{a\alpha}\pd{L}{y^\gamma_a}
-\frac{d}{dx^a}\left(\pd{L}{y^\alpha_a}\right)\right]
+\frac{d}{dx^a}\left(\sigma^\alpha\pd{L}{y^\alpha_a}\right).
\end{align*}
Let us denote by $\delta L$ the expression with components
$$
\delta L_\alpha
=\frac{d}{dx^a}\left(\pd{L}{y^\alpha_a}\right)
- Z^\gamma_{a\alpha}\pd{L}{y^\gamma_a}
-\rho^A_\alpha\pd{L}{u^A}
$$
and by $J_\sigma$ the $(r-1)$-form (along $\pi_1$)
$J_\sigma=\sigma^\alpha\pd{L}{y^\alpha_a}\omega_a$
with $\omega_a=i_{\pd{}{x^a}}\omega$. Then we have that 
$$
0=\frac{d}{ds}\CMcal{S}(\Phi_s)\at{s=0}
=\int_N(\delta L_\alpha\,\sigma^\alpha)\omega+\int_Nd(J_\sigma\circ\Phi).
$$
Since $\sigma$ has compact support the second term vanishes by Stokes theorem. Moreover, since the section $\sigma$ is arbitrary, by the fundamental theorem of the Calculus of Variations, we get $\delta L=0$, which are the Euler-Lagrange equations. Notice that the first two equations in the above statement are but the morphism conditions.
\end{proof}

\begin{remark} 
We can consider a more general problem and we can look for the critical points of the action $\int_NL(\Phi)\,\omega$ where $\Phi$ is restricted to be an admissible map, i.e. with the constraints $u^A_{,a}=\rho^A_a+\rho^A_\alpha y^\alpha_a$. If we apply formally the Lagrange multiplier trick we get the first and the third of the above equations for $\pd{L}{y^\alpha_a}=\rho^A_\alpha p^a_A$, where $p^A_a$ the Lagrange multipliers. Moreover the integrability conditions for the admissibility equations are $\rho^A_\alpha\mathcal{M}^\alpha_{ab}=0$, where 
$$
\mathcal{M}^\alpha_{ab}=\pd{y^\alpha_a}{x^b}-\pd{y^\alpha_b}{x^a}
 +C^\alpha_{b\gamma}y^\gamma_a-C^\alpha_{a\gamma}y^\gamma_b
 +C^\alpha_{\beta\gamma}y^\beta_by^\gamma_a
 +C^\alpha_{ab}
$$
i.e. 
$\mathcal{M}^\alpha_{ab}=0$ are the morphism condition. Therefore, if $\rho|_K$ is injective we recover the three field equations. Nevertheless, even in that case, we cannot ensure that the Lagrange multiplier method captures all the solutions.
\end{remark}

\subsection*{Noether's theorem}
It is well-known that Noether's theorem is a consequence of the existence of a variational description of a problem. In the standard case~\cite{Barna-L}, when the Lagrangian is invariant by the first jet prolongation of a vertical vectorfield $Z$ then the Noether current $J=i_{Z^{(1)}}\Theta_L$ is a conserved current, that is, its pullback by any solution of the Euler-Lagrange equations is a closed form on the base manifold $N$, and therefore its integral over any closed $(r-1)$-dimensional submanifold vanishes. We will show that a similar statement can be obtained for a field theory over Lie algebroids.

Consider a section $\sigma$ of $E$ vertical over $F=TN$, and the associated vector field $\XsigmaJ$. We can associate to $\sigma$ a current $J_\sigma$ (i.e. a $(r-1)$-form along~$\prJN$) as in the standard case by means of 
$$
J_\sigma=\pd{L}{y^\alpha_a}\sigma^\alpha\omega_a,
$$
which can be defined intrinsically via vertical lifting.

\begin{definition}
We will say that the Lagrangian density $L$ is \emph{invariant} under a $\prEF$-vertical section $\sigma\in\sec{E}$ if $\mathcal{L}_{\XsigmaJ}L=0$. An $(r-1)$-form $\lambda$ along $\prJN$ is said to be a \emph{conserved current} if $\lambda\circ\Phi$ is a closed form on $N$ for any solution $\Phi$ of the Euler-Lagrange equations.
\end{definition}

It follows from this definition that $L$ is invariant under $\sigma$ if and only if it is invariant under $\Jprol{\Psi_s}$, where $\Psi_s$ is the flow defined by the section $\sigma$.

\begin{theorem}
Let $\sigma\in\sec{E}$ be a $\prEF$-vertical section. If the Lagrangian  is invariant under the section $\sigma$ then $J_\sigma$ is a conserved current.
\end{theorem}
\begin{proof}
Indeed, following the steps in the derivation of the Euler Lagrange equations we have that
$
(\mathcal{L}_{\XsigmaJ}L)\omega=\delta L_\alpha\sigma^\alpha\omega 
+d(J_\sigma\circ\Phi)
$.
Therefore, if $\Phi$ is a solution of the Euler-Lagrange equations, the first term vanishes and, since $\mathcal{L}_{\XsigmaJ}L=0$, we have that $d(J_\sigma\circ\Phi)=0$.
\end{proof}

\section{Examples}
\label{examples}

\subsection*{Standard case} In the standard case, we consider a bundle $\map{\nu}{M}{N}$, the standard Lie algebroids $F=TN$ and $E=TM$ and the tangent map $\map{\pi=T\nu}{TM}{TN}$. Then we have that $\Jpi=J^1\nu$. When we choose coordinate basis of vectorfields (i.e. of sections of $TN$ and $TM$) we recover the equations for the standard first order field theory. 

Moreover, if we consider a different basis, what we get is the equations for a first order field theory but  written in pseudo-coordinates. In particular, one can take an Ehresmann connection on the bundle $\map{\nu}{M}{N}$ and use an adapted local basis 
$$
\bar{e}_i=\pd{}{x^i},
\qquand
e_i=\pd{}{x^i}+\Gamma^A_i\pd{}{u^A},
\qquad
e_A=\pd{}{u^A}.
$$
In this case, the greek indices and the latin capital indices coincide. We have the brackets 
$$
[e_i,e_j]=-R^A_{ij}e_A,
\qquad
[e_i,e_B]=\Gamma^A_{iB}e_A
\qquand
[e_A,e_B]=0,
$$
where we have written $\Gamma^B_{iA} = \partial{\Gamma^B_i}/\partial{u^A}$ and where $R^A_{ij}$ is the curvature tensor of the nonlinear connection we have chosen. The components of the anchor are $\rho^i_j=\delta^i_j$, $\rho^A_i=\Gamma^A_i$ and $\rho^A_B=\delta^A_B$ so that the Euler-Lagrange equations are 
\begin{align*}
&\pd{u^A}{x^i}=\Gamma^A_i+y^A_i\\
&\pd{y^A_{i}}{x^j}-\pd{y^A_{j}}{x^i}
+\Gamma^A_{jB}y^B_i-\Gamma^A_{iB}y^B_j
=R^A_{ij}\\
&\frac{d}{dx^i}\left(\pd{L}{y^A_i}\right)
-\Gamma^B_{iA}\pd{L}{y^B_i}
=\pd{L}{u^A}.
\end{align*} 
See~\cite{CaCrIb,Barna-L} for an alternative derivation of this equations.

\subsection*{Time-dependent Mechanics} In~\cite{SaMeMa1,SaMeMa2} we developed a theory of Lagrangian Mechanics for time dependent systems defined on Lie algebroids, where the base manifold is fibered over the real line $\Real$. Since time-dependent mechanics is but a 1-dimensional field theory, our results must reduce to that. 

The morphism condition is just the admissibility condition so that, if we write $x^0\equiv t$ and $y^\alpha_0\equiv y^\alpha$, the Euler-Lagrange equations are
$$
\begin{aligned}
&\frac{du^A}{dt}=\rho^A_0+\rho^A_\alpha y^\alpha\\
&\frac{d}{dt}\left(\pd{L}{y^\alpha}\right)
=\pd{L}{y^\gamma}(C^\gamma_{0\alpha}+C^\gamma_{\beta\alpha}y^\beta)
+\pd{L}{u^A}\rho^A_\alpha,
\end{aligned}
$$
in full agreement with~\cite{MaMeSa}. 

\subsection*{Chern-Simons theory} 

We consider a Lie algebra $\mathfrak{g}$ with an $\mathrm{ad}$-invariant metric $\boldsymbol{k}$. We choose a basis $\{\epsilon_\alpha\}$ of $\mathfrak{g}$ and the structure constants $C^\alpha_{\beta\gamma}$ in that basis. If $k_{\alpha\beta}$ are the components of the metric $\boldsymbol{k}$, then the symbols $C_{\alpha\beta\gamma}=k_{\alpha\mu}C^\mu_{\beta\gamma}$ are skewsymmetric in all indices. We consider a 3-dimensional manifold $N$ and the Lie algebroid $E=TN\times\mathfrak{g}\to N$ with the projection onto the first factor. A basis for sections of $E$ is given by $e_\alpha(n)=(n,\epsilon_\alpha)$. The maps $\pi$ and $\nu$ are 
$\pi(v,\xi)=v$ and $\nu=\id_N$. 

If a map  $\Phi$ from $TN$ to $E$ is a section of $\pi$, then $\Phi$ is of the form $\Phi(v)=(v,A^\alpha(v)\epsilon_\alpha)$ for some 1-forms $A^\alpha$ on $N$. In other words $\Phi\pb e^\alpha=A^\alpha$. The Lagrangian density for Chern-Simons theory is
$$
L\,dx^1\wedge dx^2\wedge dx^3 =\frac{1}{3!}\,C_{\alpha\beta\gamma}\,A^\alpha\wedge A^\beta\wedge A^\gamma.
$$
With the notation of this paper, the coordinates $y^\alpha_a$ are the components of $A^\alpha$, that is, $A^\alpha=y^\alpha_adx^a$. Therefore the Lagrangian is  
$L=C_{\alpha\beta\gamma}y^\alpha_1y^\alpha_2y^\alpha_3$.

There is no admissibility condition in this case, since there are no coordinates~$u^A$. The morphism conditions are $\p{y}^\alpha_{i|j}-\p{y}^\alpha_{j|i}-C^\alpha_{\beta\gamma}y^\beta_jy^\gamma_i=0$, and can be written conveniently in terms of the 1-forms $A^\alpha$ as
$$
dA^\alpha+\frac{1}{2}C^\alpha_{\beta\gamma}A^\beta\wedge A^\gamma=0.
$$

For the momenta, we have that 
$$
\pd{L}{y^\alpha_1}=C_{\alpha\beta\gamma}y^\beta_2y^\gamma_3,
\qquad
\pd{L}{y^\alpha_2}=C_{\beta\alpha\gamma}y^\beta_1y^\gamma_3,
\qquand
\pd{L}{y^\alpha_1}=C_{\beta\gamma\alpha}y^\beta_1y^\gamma_2,
$$
so that the Euler-Lagrange equations reduce to
\begin{align*}
\frac{d}{dx^a}\pd{L}{y^\alpha_a}-\pd{L}{y^\gamma_a}C^\gamma_{\beta\alpha}y^\beta_a=
C_{\alpha\beta\gamma}\Bigl[
&\hphantom{{}+{}}
(y^\beta_{2|1}-y^\beta_{1|2}-C^\beta_{\mu\nu}y^\mu_1y^\nu_2)y^\gamma_3+\\
&+(y^\beta_{1|3}-y^\beta_{3|1}-C^\beta_{\mu\nu}y^\mu_3y^\nu_1)y^\gamma_2+\\
&+(y^\beta_{3|2}-y^\beta_{2|3}-C^\beta_{\mu\nu}y^\mu_2y^\nu_3)y^\gamma_1
\hphantom{{}+{}}\Bigr]=0
\end{align*}
vanish identically in view of the morphism condition.

The conventional Lagrangian density for the Chern-Simons theory is
$$
L'\omega=k_{\alpha\beta}\left(
A^\alpha\wedge dA^\beta
+\frac{1}{3}C^\beta_{\mu\nu}A^\alpha\wedge A^\mu\wedge A^\nu
\right)
$$
The difference between $L'$ and $L$ is a multiple of the morphism condition
$$
L'\omega-L\omega=
k_{\alpha\mu}A^\mu\left[
dA^\alpha+\frac{1}{2}C^\alpha_{\beta\gamma}A^\beta\wedge A^\gamma
\right].
$$
Therefore both Lagrangians coincide on the set $\Mor{\pi}$ of morphisms, which is the set where the action is defined.

\subsection*{Systems with symmetry} The case of a system with symmetry is very important in Physics. We consider a principal bundle $\map{\nu}{P}{M}$ with structure group $G$ and we set $N=M$, $F=TN$ and $E=TP/G$ (the Atiyah algebroid of $P$), with $\pi([v])=T\nu(v)$. Sections of $\pi$ are just principal connections on $P$ and a section is a morphism if and only if it is a flat connection. The kernel $K$ is just the adjoint bundle $(P\times\mathfrak{g})/G\to M$. 

An adequate choice of a local basis of sections of $F$, $K$ and $E$ is as follows. Take a coordinate basis $\bar{e_i}=\partial/\partial x^i$ for $F=TM$, take a basis $\{\epsilon_\alpha\}$ of the Lie algebra $\mathfrak{g}$ and the corresponding sections of the adjoint bundle $\{e_\alpha\}$, so that $C^\gamma_{\alpha\beta}$ are the structure constants of the Lie algebra. Finally we take sections $\{e_i\}$ of $E$ such that $e_i$ projects to $\bar{e}_i$, that is, we chose a (local) connection and $e_i$ is the horizontal lift of $\bar{e}_i$. Thus we have $[e_i,e_j]=-\Omega^\alpha_{ij}e_\alpha$ and $[e_i,e_\alpha]=0$.
In this case there are no coordinates $u^A$, and  with the above choice of basis we have that the Euler-Lagrange equations are
\begin{align*}
& \pd{y^\alpha_a}{x^b}-\pd{y^\alpha_b}{x^a}
 +C^\alpha_{\beta\gamma}y^\beta_by^\gamma_a
 =\Omega^\alpha_{ab}\\
&\frac{d}{dx^a}\left(\pd{L}{y^\alpha_a}\right)
-\pd{L}{y^\gamma_a}C^\gamma_{\beta\alpha}y^\beta_a=0.
\end{align*}

In particular, for the definition of our sections we can choose a flat connection (for instance a solution of our variational problem or just a coordinate basis). Then we have that $\Omega^\alpha_{ab}=0$ in the above equations, which are then called the covariant Euler-Poincaré equations~\cite{CaRaSh,CaGaRa}.

\section{Conclusions and outlook}
In this contribution, we have dealt with the constrained variational problem consisting in finding the critical sections of the integral of a Lagrangian function, where the fields are restricted to be morphisms from $TN$ to a Lie algebroid $E$, and we have found the Euler-Lagrange equations. The admissible variations are determined by the geometry of the problem and are not prescribed in an ad hoc manner.

Particular cases of our Euler-Lagrange equations are the Euler-Poincaré equations for a system defined on the bundle of connections of a principal bundle, the Lagrange-Poincaré equations for systems defined on a bundle of homogeneous spaces and the Lagrange-Poincaré for systems defined in semidirect products. 

In the standard case, when $E=TM$, the existence of a multisymplectic form is fundamental in the development of the theory. Preliminary calculations show that in the case of arbitrary Lie algebroids one can also find a Cartan form~\cite{Garcia,GoSt} and a multisymplectic form, and the field equations are obtained via a multisymplectic equation.

Another interesting problem under development is the determination of those Lagrangians which give rise to null Euler-Lagrange equations, as it was the case of the  Chern-Simons theory. Finally, it would be interesting to find a variational principle for the case of a general Lie algebroid $F\not=TN$.

\end{document}